\documentclass[11pt]{article}

\usepackage{amsmath}
\usepackage{amsmath,amsfonts,amssymb,amsthm,graphics,amscd}
\usepackage[lmargin=3cm, rmargin=3cm,tmargin=2.5cm,bmargin=2.5cm]{geometry}

\usepackage{color}

\newcommand{\CC}{{\mathbb C}}

\newcommand{\bdf}{\begin{definition}}
\newcommand{\edf}{\end{definition}}
\newcommand{\bprop}{\begin{proposition}}
\newcommand{\eprop}{\end{proposition}}
\newcommand{\bcor}{\begin{corollary}}
\newcommand{\ecor}{\end{corollary}}
\newcommand{\bet}{\begin{theorem}}
\newcommand{\eet}{\end{theorem}}
\newcommand{\blm}{\begin{lemma}}
\newcommand{\elm}{\end{lemma}}
\newcommand{\bp}{\begin{proof}}
\newcommand{\ep}{\end{proof}}
\newcommand{\bex}{\begin{example}\rm}
\newcommand{\eex}{\end{example}}
\newcommand{\bexs}{\begin{examples}\rm}
\newcommand{\eexs}{\end{examples}}
\newcommand{\bremark}{\begin{remark}\rm}
\newcommand{\eremark}{\end{remark}}

\newcommand{\norm}[1]{ \| #1 \| }

\newtheorem{theorem}{Theorem}[section]
\newtheorem{lemma}[theorem]{Lemma}
\newtheorem{corollary}[theorem]{Corollary}

\newtheorem{example}[theorem]{Example}
\newtheorem{definition}[theorem]{Definition}
\newtheorem{proposition}[theorem]{Proposition}
\newtheorem{remark}[theorem]{Remark}

\newtheorem{examples}[theorem]{Examples}

\def\CC{\mathbb{C}}

\numberwithin{equation}{section}

\def\CC{{\mathbb C}}

\begin{document}

\title {On   Classes of Fredholm Type  Operators}

\author { Alaa Hamdan and  Mohammed
Berkani\ }

\date{}

\maketitle \setcounter{page}{1}

\begin{abstract}
Given  an  idempotent $p$ in a Banach algebra and following the study in \cite{P50} of p-invertibility, we
consider here left p-invertibility, right p-invertibility  and p-invertibility in the Calkin Algebra $\mathcal{C}(X),$
where $X$ is a Banach space. Then we  define and study  left and right   generalized Drazin invertibility
and we characterize left and right Drazin invertible elements in the Calkin  algebra. Globally, this leads to define and  characterize the classes of P-Fredholm, pseudo B-Fredholm and weak B-Fredholm operators.

 \end{abstract}

\renewcommand{\thefootnote}{}
\footnotetext{\hspace{-7pt} {\em 2020 Mathematics Subject
Classification\/}: 47A53,  16U90, 16U40.
\baselineskip=18pt\newline\indent {\em Key words and phrases\/}:
Fredholm operators, Calkin algebra,  generalized Drazin invertibility,
p-invertibility.}

\section{Introduction}

Let $X$ be a Banach space,
let $L(X)$  be the Banach algebra of bounded linear operators acting on
the Banach space $X$ and let  $ T \in L(X).$
We will denote by $N(T)$  the null space of  $T$,  by $ \alpha(T)$
the nullity of $T$,  by $ R(T)$ the range of $T$ and by $\beta(T) $
its defect. If the range $ R(T)$ of $T$ is closed  and  $ \alpha(T) < \infty $
 (resp.  $\beta(T)< \infty $ ),
\noindent  then $T$ is called an upper semi-Fredholm (resp. a
lower semi-Fredholm) operator. A semi-Fredholm operator is an
upper or a lower semi-Fredholm operator. If both of $ \alpha(T) $
and $\beta(T) $ are finite then $T $ is called a Fredholm operator
and the index of $T$  is defined by   $ ind(T) = \alpha(T) -
\beta(T). $  The notations $\Phi_+(X),  \Phi_-(X)$ and  $\Phi(X)$ will design  respectively
the set of upper semi-Fredholm, lower semi-Fredholm and Fredholm operators. 

For $T\in L(X),$ the Fredholm  spectrum of $T$ is defined by:  $$ \sigma_F(T) = \{\lambda \in \CC: T-\lambda I \text { is\ not\ a\ Fredholm\ operator}  \}.$$
Define also  the sets:

$\Phi_l(X) = \{T \in \Phi_+ (X) \mid $ there exists a bounded projection of $X$ onto $R(T)$\},\, 

and

$\Phi_r (X) = \{T \in \Phi_-(X)\mid $ there exists a bounded projection of $X $ onto $N(T) \}$.\\

\noindent Recall that the Calkin algebra over $X$ is the quotient algebra
$\mathcal{C}(X)=L(X)/K(X)$, where  $K(X)$ is the closed ideal of
compact operators on $X$.  Let $G_r$ and $G_l$   be the right and left, respectively, invertible elements
  of $ \mathcal{C}(X)$. From \cite[Theorem 4.3.2]{CPY}    and  \cite[Theorem 4.3.3]{CPY} ,
   it follows that $ \Phi_l (X)= \Pi^{-1}(G_l)$  and  $ \Phi_r (X)= \Pi^{-1}(G_r),$
     where $ \Pi: L(X) \rightarrow  \mathcal{C}(X)$ is the natural projection. We observe that $\Phi(X)= \Phi_l (X)\cap \Phi_r(X).$

  \bdf The elements of  $ \Phi_l (X)$ and   $\Phi_r (X)$ will be  called respectively left
    semi-Fredholm operators and   right  semi-Fredholm operators .

  \edf

In 1958, in his paper \cite{DR},  M.P. Drazin  extended the
concept of invertibility in  rings and semigroups and introduced a
new kind of inverse, known now as the  Drazin inverse.

\bdf An element  $a$  of  a semigroup $ \mathcal{S}$ is called
Drazin invertible if there exists an element $b \in \mathcal{S}$
written $ b= a^d$  and called the Drazin inverse of $a,$
satisfying the following equations: $$ ab = ba, b=ab^2, a^k
 = a^{k+1}b,$$
 for some nonnegative integer $k.$ The least nonnegative
integer $k$ for which these equations holds is the Drazin index $
i(a)$ of $a.$

\edf

\noindent It follows from \cite{DR} that a Drazin invertible
element in a semigroup has a unique Drazin inverse.\\

\noindent In 1996, in  \cite[Definition 2.3]{Koliha}, J.J. Koliha extended
the notion of  Drazin invertibility.

\bdf An element  $a$  of a Banach algebra $A$ will be said to be
generalized Drazin invertible if there exists $b \in A $ such that
$ bab=b, ab=ba$ and  $ aba-a$ is a quasinilpotent element in $A$.
\edf

Koliha \cite{Koliha} proved that  $a\in A$ is generalized Drazin
invertible if and only if  there exists $ \epsilon >0,$ such that
for all $\lambda $ such that $ 0 < \mid \lambda \mid <  \epsilon,
$ the element  $a-\lambda e $ is invertible and he proved in \cite[Theorem 4.2]{Koliha},
 that a generalized Drazin
invertible element  has a unique generalized Drazin
inverse. He also proved that an element  $a\in A$ is generalized Drazin
invertible if and only if  there exists an idempotent $p \in A$  commuting with $a,$
such that $ a+p$ is invertible in $A$ and $ap$ is quasinilpotent. \\

Let $A$ be a ring with a unit and let $p$ be an idempotent in $A$.  Recall that the  commutant $C_p$ of $p$  is the subring of $A$ defined  by $C_p=
\{ x \in A \mid xp=px \}.$
   In  \cite[Definition 2.2]{P50},
the concepts of left $p-$invertibility, right $p-$invertibility and $p-$invertibility where defined  as follows.

\bdf Let $ a \in A.$   We will say that:

1- $a$ is left  $p$-invertible if $ ap=pa$ and  $ a+p$ is left
invertible in $C_p.$

2- $a$ is right  $p$-invertible if $ ap=pa$ and $ a+p$ is right
invertible in $C_p.$

 3- $a$ is
$p$-invertible if $ ap=pa$    and $ a+p$ is invertible.

\edf

Moreover in \cite[Definition 2.11]{P50}, left and right Drazin invertibility where defined as follows:

\bdf \label{Drazin $p$-invertibility}  We will say that an element
$a  \in A $ is left Drazin invertible (respectively right Drazin
invertible) if there exists an idempotent $p \in A$ such that $a$
is left $p$-invertible (respectively right $p$-invertible) and
$ap$ is nilpotent.

\edf

For  $T \in L(X),$  we will say that  a subspace $M$
 of $X$  is {\em $T$-invariant} if $T(M) \subset M.$ We define
$T_{\mid M}:M \to M$ as $T_{\mid M}(x)=T(x), \, x \in M$.  If $M$
and $N$ are two closed $T$-invariant subspaces of $X$ such that
$X=M \oplus N$, we say that $T$ is {\em completely reduced} by the
pair $(M,N)$ and it is denoted by $(M,N) \in Red(T)$. In this case
we write $T=T_{\mid M} \oplus T_{\mid N}$ and we say that $T$ is a
{\em direct sum} of $T_{\mid M}$ and $T_{\mid N}$.

It is said that $T \in L(X)$ admits a  generalized
Kato decomposition, abbreviated as GKD,  if there exists $(M, N)
\in Red(T)$ such that $ T_{\mid M}$ is Kato and $T_{\mid N}$ is
quasinilpotent. Recall that an operator $T \in L(X)$ is {\em Kato}
if $R(T)$ is closed and $Ker(T) \subset R(T^n)$ for every $ n \in
\mathbb{N}$.
\bremark For $T \in L(X),$ to say that a pair $(M,N)$ of closed subspaces of $X$ is  in $Red(T),$
means simply that there exists an idempotent $P \in L(X)$ commuting with $T.$ Indeed if $(M,N) \in Red(T),$
let $P$ be the projection on $M$ in parallel to $N.$ Then since $(M,N) \in Red(T),$ we see that  $P$ commutes  with $T.$
Conversely, given an idempotent $P \in L(X)$ commuting with $T,$ if we set $M= P(X)$ and $N= (I-P)X,$ it is clear that $(M,N) \in Red(T).$

\eremark

For $T \in L(X)$ and a nonnegative integer $n,$ define $T_{[n]}$ to
be the restriction of $T$ to $R(T^n)$ viewed as a map from $R(T^n)$
into $R(T^n)$ ( in particular $T_{[0]}=T).$ If for some integer $n$
the range space $R(T^n)$ is closed and $T_{[n]}$ is an upper (resp.
a lower) semi-Fredholm operator, then $T$ is called an  upper (resp.
a lower) semi-B-Fredholm operator.
Moreover, if $T_{[n]}$ is a Fredholm operator, then $T$ is called a
B-Fredholm operator, see \cite{P7}. From    \cite[Theorem 2.7 ]{P7}, we know  that   $T \in L(X)$  is  a B-Fredholm operator if
 there exists $(M, N) \in Red(T)$ such that $T_{\mid M}$ is a nilpotent   operator and $T_{\mid N}$ is a Fredholm operator.  The notations $ \Phi_B(X), \Phi_{B^+}(X)$ and  $\Phi_{B^-}(X)$ will designate  respectively the set of upper semi-B-Fredholm, lower semi-B-Fredholm and B-Fredholm operators. We have $ \Phi_B(X)= \Phi_{B^+}(X) \cap \Phi_{B^-}(X). $

 The results of this paper can be summarized as follows. First let us recall the following important result, which we will use throughout this paper. If $p$ is an idempotent element of $\mathcal{C}(X),$ then from
   \cite [Lemma 1]{BAR}    there exists an idempotent $P \in L(X)$ such that $ \Pi(P)= p.$ So letting $p=\Pi(P)$ being  an idempotent  in the Calkin algebra and    using the results of \cite{P50},   we characterize in the second section   algebraically
 $p$-invertible (respectively right $p$-invertible and  left $p$-invertible) elements in the Calkin algebra
 $\mathcal{C}(X).$ This analysis leads to the introduction of the class $\Phi_P(X)$ of $P$-Fredholm ( respectively the class $\Phi_{lP}(X)$ of left semi-$P$-Fredholm and the class $\Phi_{rP}(X)$ of right  semi-$P$-Fredholm)  operators in a similar way as  the corresponding classes of Fredholm, left semi-Fredholm and right semi-Fredholm operators.
We will see that $\Phi_P(X)= \Phi_{lP}(X)\cap \Phi_{rP}(X).$

 In the third section, we introduce the notions  of left and right generalized invertibility
 in a Banach algebra and we   characterize left and right generalized Drazin invertible elements in the Calkin algebra. In particular  we prove that an element of $\mathcal{C}(X)$ is generalized Drazin invertible if and only if it is left and right generalized Drazin invertible. This leads us  to introduce the classes $\Phi_{l\mathcal{P}B}(X) $ and $\Phi_{r\mathcal{P}B}(X) $ of left and right  pseudo semi-B-Fredholm operators, completing in this way the study of the class  $\Phi_{\mathcal{PB}}(X) $ of pseudo B-Fredholm operators inaugurated  in \cite{P45}. We will prove  $\Phi_{\mathcal{P}B}(X)= \Phi_{l\mathcal{P}B}(X)\cap \Phi_{r\mathcal{P}B}(X).$

 Similarly in the fourth section, after characterizing left and right  Drazin invertible elements in the Calkin algebra, we prove that an element of $\mathcal{C}(X)$ is  Drazin invertible if and only if it is left and right  Drazin invertible.  From \cite[Theorem 3.4]{P10}, we know that  $T \in L(X)$ is a B-Fredholm operator
if and only if $\pi(T)$ is Drazin invertible in the algebra $
L(X)/F_0(X), $ where $F_0(X)$ is the ideal of finite rank operators in $L(X),$ and $\pi: L(X) \rightarrow L(X)/F_0(X)$ is the natural projection. Thus to differentiate the class of operators studied in this section from B-Fredholm operators,  we name them  ''weak B-Fredholm operators''. We will prove that the class of weak B-Fredholm operators contains strictly the class of B-Fredholm operators.

This leads us  to introduce the classes $\Phi_{l\mathcal{W}B}(X),$  $\Phi_{r\mathcal{W}B}(X) $  of  left and right  weak semi-B-Fredholm operators, completing in this way the study of the class  $\Phi_{\mathcal{W}B}(X) $ of weak B-Fredholm operators inaugurated  in \cite{P45}. We will prove that  $\Phi_{\mathcal{W}B}(X)= \Phi_{l\mathcal{W}B}(X)\cap \Phi_{r\mathcal{W}B}(X).$

\section{ On P-Fredholm Operators}

\bdf Let $ T \in L(X)$   and let $p= \Pi(P)$ be an idempotent element of $\mathcal{C}(X).$
We will say that:

1- $T$ is a left semi-$P$-Fredholm operator if $ \Pi(T)  $ is left $p$-invertible in $\mathcal{C}(X).$

2- $T$ is a right  semi-$P$-Fredholm operator if $ \Pi(T)$ is right $p$-invertible in $\mathcal{C}(X).$

 3- $T$ is a $P$-Fredholm operator if $ \Pi(T)$ is  $p$-invertible in $\mathcal{C}(X)$.

\edf

Let $p$ be an idempotent element of $\mathcal{C}(X),$ we  define the sets  $G_{rp}$ and $G_{lp}$  to  be the right and left, respectively, p-invertible elements  of $ \mathcal{C}(X).$    Then  $ \Phi_{lP} (X)= \Pi^{-1}(G_{lp})$  and  $ \Phi_{rP} (X)= \Pi^{-1}(G_{rp}).$     If  the idempotent p=0,  then $ \Phi_{l0} (X)=  \Phi_l (X),  \Phi_{r0} (X)=  \Phi_r (X),       G_{l0}=  G_l$  and   $G_{r0}=  G_r .$
So we retrieve the previous  relations $ \Phi_l (X)= \Pi^{-1}(G_l)$  and  $ \Phi_r (X)= \Pi^{-1}(G_r).$\\

The following results characterizes algebraically left  semi-$P$-Fredholm,   right  semi-$P$-Fredholm  and $P$-Fredholm operators. Since these characterizations are obtained straightforwardly by using the corresponding results obtained in  \cite{P50}, we give them without proofs, referring the reader to  \cite{P50}. \\

Using \cite[Theorem 2.5]{P50},  we obtain  the following characterization of  left  semi-$p$-Fredholm operators.

\bet Let $ T \in L(X).$  Then there exist an idempotent $p=\Pi(P) \in \mathcal{C}(X)$ such that  $T$ is a left  semi-P-Fredholm operator if and
only if there exists an element  $S \in L(X)$ such that $STS-S, S^2T-S, TST-ST^2$ are compact, there exists $ V \in L(X)$ such that the commutator $[V, ST]$  and  $ V(I+T-ST) -I$ are compact. In this case  $ p= \Pi( I-ST).$

\eet

Using \cite[Theorem 2.6]{P50},  we obtain  the following characterization of  right  semi-$P$-Fredholm operators..

\bet Let $ T \in L(X).$  Then  there exist an idempotent $p=\Pi(P) \in \mathcal{C}(X)$ such that $T$ is a right   semi-P-Fredholm operator  if and
only if there exists an element  $S \in L(X)$ such that $STS-S, TS^2-S, TST- T^2S $ are compact, there exists $ V \in L(X)$ such that the commutator  $[V, TS]$  and  $ (I+T-TS)V -I$ are compact. In this case  $ p= \Pi( I-TS).$

\eet

Using \cite[Theorem 2.7]{P50}, we obtain  the following characterization of $P$-Fredholm operators..

\bet Let $ T \in L(X).$  Then there exist an idempotent $p=\Pi(P) \in \mathcal{C}(X)$ such that $T$ is a $P$-Fredholm operator  if and
only if there exists an element  $S \in L(X)$ such that the commutator $[T,S]$ and $STS-S$  are compact and
that $  I+T- ST$ is  a Fredholm operator. In this case  $ p= \Pi( I-ST).$

\eet

\bremark Let $p=\Pi(P) \in \mathcal{C}(X)$ be an idempotent. It is clear that  $ T\in L(X)$ is $P$-Fredholm   if and only if $ T\in L(X)$ is left and right semi-$P$-Fredholm. So we have  $\Phi_P(X)= \Phi_{lP}(X)\cap \Phi_{rP}(X).$
\eremark

\section{ On Pseudo B-Fredholm Operators}

\bdf

\begin{enumerate}
  \item $T \in L(X)$  is called a Riesz- left semi-Fredholm operator if
 there exists $(M, N) \in Red(T)$ such that $
T_{\mid M}$ is a Riesz   operator and $T_{\mid N} \in \Phi_l(N) .$
  \item  $T \in L(X)$  is called a Riesz- right  semi-Fredholm operator if
 there exists $(M, N) \in Red(T)$ such that $
T_{\mid M}$ is a Riesz  operator and $T_{\mid N} \in \Phi_r(N) .$

  \item  $T \in L(X)$  is called a Riesz- Fredholm operator if
 there exists $(M, N) \in Red(T)$ such that $
T_{\mid M}$ is a Riesz  operator and $T_{\mid N}$ is a Fredholm operator.
\end{enumerate}
\edf

\bdf \label{Drazin $p$-invertibility}  Let $A$ be a Banach algebra. We will say that an element
$a  \in A $ is:

\begin{enumerate}
\item left generalized Drazin invertible  if there exists an idempotent $p \in A$ such that $a$
is left $p$-invertible in $A$ and $ap$ is quasinilpotent.

\item   right generalized  Drazin
invertible if there exists an idempotent $p \in A$ such that $a$
is  right $p$-invertible in $A$ and $ap$ is quasinilpotent.

\item  generalized Drazin invertible  if there exists an idempotent $p \in A$  such that $a$
is  $p$-invertible in $A$  and
$ap$ is quasinilpotent.

\end{enumerate}

\edf

\bdf \label{pseudo} Let $ T \in L(X).$    We will say that:

\begin{enumerate}
  \item $T$ is a left pseudo   semi-B-Fredholm operator if $ \Pi(T)  $ is left generalized Drazin  invertible in $\mathcal{C}(X).$
  \item $T$ is a right pseudo semi-B-Fredholm operator if $ \Pi(T)$ is right generalized  Drazin invertible in $\mathcal{C}(X).$
  \item $T$ is a  pseudo semi-B-Fredholm operator if it is right  or left pseudo semi-B-Fredholm.
  \item $T$ is a pseudo B-Fredholm operator if $ \Pi(T)$ is  generalized Drazin invertible in $\mathcal{C}(X).$
\end{enumerate}

\edf

\begin{remark}

\begin{enumerate}

\item The class of  pseudo B-Fredholm  was defined in  \cite[Definition 2.4]{P45}. It involves the class
of B-Fredholm operators defined in \cite{P7}.

 \item In \cite{BO}, in 2015, the author studied  generalized Drazin invertible elements under Banach algebra homomorphisms.  In the case of the  Calkin Algebra,  he proved in \cite[Theorem 6.1]{BO} that an operator whose image is generalized Drazin invertible in the Calkin algebra is a Riesz-Fredholm operator. Later in 2019, in \cite[Theorem 2.10]{P45},  the authors proved the same result following using a shortened proof. 

 Here in this paper,  the concept of p-invertibility enable us to give left and right version of \cite[Theorem 6.1]{BO}.

\end{enumerate}

\end{remark}

\bet \label{left-pseudo} Let $T \in L(X).$ Then the following properties
are equivalent:

\begin{enumerate}

\item $T$ is a left pseudo semi-B-Fredholm operator.

\item  $T$ is a compact perturbation of a Riesz- left semi-Fredholm operator.

\end{enumerate}

\eet
\bp  $1)\Rightarrow 2)$  Assume that  $T$ is a left pseudo semi-B-Fredholm operator, then
$\Pi(T)$  is left generalized Drazin invertible in  $\mathcal{C}(X).$ So there exist an
idempotent $p =\Pi(P) \in  \mathcal{C}(X)$ such that:

\begin{itemize}
  \item  $p \Pi(T) = \Pi(T) p$
  \item $p \Pi(T)$ is quasinilpotent in $\mathcal{C}(X)$
  \item  There exists  $S \in L(X)$  such that $p \Pi(S) = \Pi(S) p$ and  $\Pi(S) ( \Pi(T) + p)) = \Pi(I)$

\end{itemize}

 Let $ X= X_1 \oplus X_2$ be
the decomposition of $X$ associated to $P,$ that's  $ X_1= R(P)$ and $ X_2= N(P).$
 Since $\Pi(T)$ and $\Pi(P)$ commutes, we have that $\Pi(PTP)=\Pi(TP)$ and $\Pi((I-P)T(I-P))=\Pi(T(I-P))$. It follows that $PTP$
 is a Riesz operator, $TP=PTP+K_1$, $T(I-P)=(I-P)T(I-P)+K_2$, where $K_1,K_2\in K(X)$, and so:
  \begin{eqnarray*}
    T &=& TP+ T(I-P)=PTP+(I-P)T(I-P)+K,
  \end{eqnarray*}
  \noindent where $K=K_1+K_2\in K(X)$.
    We have that $(X_1, X_2)\in Red(PTP)$ and  $(X_1,X_2)\in Red((I-P)T(I-P))$, hence:
  \begin{equation*}
    PTP = (PTP)_{\mid X_1}\oplus (PTP)_{\mid X_2}=(PTP)_{\mid X_1}\oplus 0,
  \end{equation*}
  and

 \begin{equation*}
      (I-P)T(I-P) =((I-P)T(I-P)) _{\mid X_1}\oplus ((I-P)T(I-P))_{\mid X_2}=0\oplus ((I-P)T(I-P))_{\mid X_2},
  \end{equation*}

  Therefore,

\begin{equation}\label{first-1}
    T=(PTP)_{\mid X_1}\oplus ((I-P)T(I-P))_{\mid X_2}+K.
  \end{equation}

  It's easily seen that $(PTP)_{\mid R(P)}$ is Riesz operator. Moreover  we show that $((I-P)T(I-P))_{\mid R(I-P)}$ is a left semi-Fredholm operator.
Since $\Pi(S) ( \Pi(T) + p) = \Pi(I),$ we have: $$ (I-P)S(I-P) (I-P) (T+P)(I-P)= I-P  +F_1,$$
where   $F_1$ is a  compact  operator such that $ (X_1, X_2) \in Red(F_1).$
As $I-P$ is the identity on $L(X_2)$,  it follows   that
$((I-P)T(I-P))_{\mid X_2}=((I-P)(T+P)(I-P))_{\mid X_2}$ is a left semi-Fredholm operator. According to (\ref{first-1}), we see that  $T$ is a compact perturbation of
a Riesz-left  semi-Fredholm operator.

$ 2)\Rightarrow 1)$ Conversely assume that $ X= X_1 \oplus X_2$  is the direct sum of the two closed subspaces $X_1, X_2$  of $X$ such that  $T=T_1\oplus T_2+ K $ where $T_1 \in L(X_1) $ is a Riesz operator and $T_2 \in L(X_2)$ is a left semi-Fredholm operator and $K\in K(X)$. Let $P$ be the projection of $X$ on $X_1$ in parallel to $X_2.$
Then $ PT= (T_1 \oplus 0) +PK$ and $TP= (T_1 \oplus 0) +KP.$  Thus $ PT-TP= PK-KP$ and $PK-KP$  is compact. Let $p= \Pi(P),$  then   $p$ commutes with $ \Pi(T).$
Let also $I_1= {P_|}_{ X_1}$  and $I_2= {(I-P)_|}_{X_2}.$
As $T_2$ is a left  semi-Fredholm operator in $L(X_2),$  there exists $ S_2 \in  L(X_2)$  such that $ S_2 T_2- I_2$ is compact and $I_1 \oplus S_2$ commutes with $P$   because
$(I_1  \oplus S_2)P=P(I_1 \oplus S_2)=I_1 \oplus 0.$

We have $ T_1\oplus T_2 +P = ( T_1  \oplus 0) +  P  + ( 0 \oplus T_2) = P[ (T_1 \oplus 0)+ I)] P  + (I-P) ( I_1 \oplus  T_2 ) (I-P).$  As $T_1$ is a Riesz operator, then $T_1 \oplus 0$ is also a Riesz operator and $\Pi( (T_1 \oplus 0)+ I))=  \Pi( (T_1 \oplus 0)) + \Pi(I)$ is invertible in $\mathcal{C}(X).$ Let $ \Pi(S_1) $ be its inverse, where $ S_1 \in L(X).$
We observe that $\Pi( I_1\oplus  T_2) $ is left invertible in  $\mathcal{C}(X)$ having  $ \Pi( I_1 \oplus S_2)$ as a left inverse.Then:

$$ \Pi((PS_1P + (I-P)(I_1\oplus S_2) (I-P))) \Pi ( T+P)= \Pi((PS_1P + (I-P)(I_1\oplus S_2) (I-P))) \Pi (T_1\oplus T_2 +P)$$
$$=\Pi([PS_1P + (I-P)(I_1\oplus S_2) (I-P)]) \Pi ( P[ (T_1 \oplus 0)+ I)] P  + (I-P) ( I_1 \oplus  T_2 ) (I-P))=\Pi(I).$$
It is easily seen that  $\Pi((PS_1P + (I-P)( I_1 \oplus S_2) (I-P)))$ commutes with $\Pi(P).$

\noindent Moreover we have $\Pi(PT)= \Pi( T_1 \oplus 0)$   is  quasinilpotent in $\mathcal{C}(X)$  because $ T_1 \oplus 0$ is a Riesz operator.  Thus $ T$ is a left pseudo semi-B-Fredholm operator
 \ep

Following the same method as in Theorem \ref{left-pseudo}, we obtain the following characterization of right pseudo semi-B-Fredholm operators. Due to the similarity of its proof with that of
Theorem \ref{left-pseudo}, we give it without proof.

\bet \label{right-pseudo} Let $T \in L(X).$ Then the following properties
are equivalent:

\begin{enumerate}

\item $T$ is a right pseudo semi-B-Fredholm operator.

\item  $T$ is a compact perturbation of a Riesz-right semi-Fredholm operator.

\end{enumerate}
\eet
Combining Theorem\ref{left-pseudo} and Theorem\ref{right-pseudo}, we obtain the following natural characterization of pseudo B-Fredholm operators.

\bet \label{first equivalence} Let $T \in L(X).$   Then the following properties
are equivalent:

\begin{enumerate}

\item $T$ is a pseudo B-Fredholm operator.

\item $T$  is left and  right pseudo semi-B-Fredholm operator.

\end{enumerate}

\eet

\bp  $1)\Rightarrow 2)$  Obviously if we assume that $T$ is a pseudo B-Fredholm operator,  then $T$  is left and  right pseudo semi-B-Fredholm operator.

$2)\Rightarrow 1)$ Assume now that $T$  is left and  right pseudo semi-B-Fredholm operator.
 So there exists an idempotent $p$  in $\mathcal{C}(X)$ such that $\Pi(T)$ is left $p$-invertible   and
$\Pi(T)p$ is quasinilpotent.  We have  $  \Pi(T) = p\Pi(T)p + (e-p)\Pi(T)(e-p)= x_1 + x_2,$  where $ x_1= p\Pi(T)p,$
$ x_2= (e-p)\Pi(T)(e-p)$  and $e= \Pi(I)$ is the identity element of $\mathcal{C}(X).$
We know that $p\, \mathcal{C}(X)p$ is a closed subalgebra of $\mathcal{C}(X)$ having $p$ as a unit element and
$(e-p) \mathcal{C}(X)(e-p)$ is a closed subalgebra of $\mathcal{C}(X)$ having $e-p$ as a unit element.
Now since $x_1$ is quasinilpotent in $p \, \mathcal{C}(X)p$ and $ x_2$ is left invertible in
 $(e-p)\mathcal{C}(X)(e-p),$ then if $ \lambda \in \mathbb{C}\setminus\{0\}$ and  $\mid \lambda \mid $ is small enough, we have
$ \Pi(T)- \lambda e= (x_1 - \lambda p) + (x_2 - \lambda (e-p))$ is left invertible in  $\mathcal{C}(X).$

Similarly as $\Pi(T)$  is  right generalized Drazin invertible in  $\mathcal{C}(X),$ we can prove that
 if $ \lambda \in \mathbb{C}\setminus\{0\}$ and  $\mid \lambda \mid $ is small enough, then $ \Pi(T)- \lambda e$ is right invertible in  $\mathcal{C}(X).$

Consequently if $ \lambda \in \mathbb{C}\setminus\{0\}$ and  $\mid \lambda \mid $ is small enough, then  $ \Pi(T)- \lambda e$ is invertible.
So if $ 0 \in \sigma(\Pi(T)),$  then it is an isolated point.  From \cite[Theorem 3.1]{Koliha}, it follows that $\Pi(T)$ is generalized Drazin invertible in $\mathcal{C}(X)$ and so  $T$ is a  pseudo B-Fredholm operator.
\ep
\bremark 
\begin{enumerate}
\item In the recent works \cite{AZOZ} and \cite{TAKAM},
the authors studied operators similar to  pseudo semi-B-Fredholm operators considered here. They defined them  as the direct sum
of a semi-Fredholm operator and a quasi-nilpotent one.
However there exists operators which are pseudo semi-B-Fredholm
operators in the sense of Definition \ref{pseudo}, but do not have a
decomposition  as the direct sum of a semi-Fredholm operator and a
quasi-nilpotent operator as seen by the following example.

\end{enumerate}

\begin{example}\label{infinite spectrum}  Let $T$ be a compact   operator having infinite spectrum. Since $
\Pi(T)= 0,$ then $T$ is a pseudo B-Fredholm operator   in the
sense of Definition \ref{pseudo}.  We prove that $T$ cannot be
written as the direct sum of a semi-Fredholm operator and a
quasi-nilpotent one.

Assume that  there exists a pair $(M, N) \in Red(T)$   such that $ T= T_1\oplus T_2 $ where
$T_1=T_{\mid M}$ is a quasi-nilpotent operator and $T_2=T_{\mid
N}$ is a semi-Fredholm operator.
As $T$ is compact,  $T_2$ is also  compact. Then  necessarily $N$ is finite dimensional and then $\sigma(T_2)$  is  finite.
Since $\sigma(T)=\sigma(T_1)\cup\sigma(T_2),$ it follows that $\sigma(T)$ is finite and this is  a contradiction.
\end{example}

\eremark

\section {On Weak B-Fredholm Operators }

\bdf

\begin{enumerate}

  \item $T \in L(X)$  is called a power compact-left semi-Fredholm operator if
 there exists $(M, N) \in Red(T)$ such that $
T_{\mid M}$ is a power compact   operator and $T_{\mid N} \in \Phi_l(N) .$

  \item $T \in L(X)$  is called  a power compact-right  semi-Fredholm operator if
 there exists $(M, N) \in Red(T)$ such that $
T_{\mid M}$ is a power compact  operator and $T_{\mid N} \in \Phi_r(N) .$

\end{enumerate}

\edf

\bdf \label{Drazin $p$-invertibility} \cite[Definition 2.11]{P50}  Let $A$ be a Banach algebra. We will say that an element
$a  \in A $ is  left  Drazin invertible (respectively right generalized  Drazin
invertible) if there exists an idempotent $p \in A$ such that $a$
is left $p$-invertible (respectively right $p$-invertible) and
$ap$ is nilpotent.

\edf

\bdf \label{Weak} Let $ T \in L(X).$    We will say that:

\begin{enumerate}
  \item $T$ is a left  weak semi-B-Fredholm operator if $ \Pi(T)  $ is left  Drazin  invertible in $\mathcal{C}(X).$
  \item $T$ is a right weak semi-B-Fredholm operator if $ \Pi(T)$ is right Drazin invertible in $\mathcal{C}(X).$
  \item $T$ is a  weak semi-B-Fredholm operator if it is right  or left weak semi-B-Fredholm operator in $\mathcal{C}(X).$
  \item $T$ is a weak B-Fredholm operator if $ \Pi(T)$ is   Drazin invertible in $\mathcal{C}(X).$
\end{enumerate}

\edf

\bremark   As each Drazin invertible element in the algebra $ L(X)/ F_0(X)$ is a Drazin invertible element in the Calkin algebra, the class of B-Fredholm operators is contained in the class of weak B-Fredholm operators.

Let  $(\lambda_n)_n $ be a
sequence in $  \mathbb{C}$ such that $\lambda_n \not=0$ for all $n$ and $
\lambda_n \longrightarrow 0$ as $n \rightarrow \infty$ and let us
consider the operator $T$ defined on the Hilbert space $l^2(\bf
\mathbb{N})$ by: $$T(\xi_1,\xi_2,\xi_3,.....)=(\lambda_1\xi_1,\lambda_2\xi_2,\lambda_3\xi_3,.....).$$

Then $$
T^n(\xi_1,\xi_2,\xi_3,.....)=((\lambda_1)^n\xi_1,(\lambda_2)^n\xi_2,(\lambda_3)^n\xi_3,.....).$$
\noindent Since $(\lambda_m)^n \not= 0$ for all $m\geq 0$ and
$(\lambda_m)^n \longrightarrow 0 $ as $m \longrightarrow \infty$
for all $n\geq 0$ then $K^n \in K(X)$  and $K^n $ is not a  finite
rank  operator for all $n \geq 1$. Hence $R(K^n)$ is not closed
for all $ n \geq 1 $. Thus $K$ is not a B-Fredholm operator.  However as $ \Pi(K)=0 $ is Drazin invertible in $\mathcal{C}(X),$ $K$ is a weak B-Fredholm operator.

Thus the class of weak B-Fredholm operators $\Phi_{\mathcal{W}B}(X)$  contains strictly the class of B-Fredholm operators $\Phi_{B}(X)$, which itself contains strictly the class of Fredholm operators $\Phi(X)$ as seen in \cite{P7} and  Obviously the class of pseudo B-Fredholm operators $\Phi_{\mathcal{P}B}(X)$ contains the class of weak B-Fredholm operators $\Phi_{\mathcal{W}B}(X)$ .

Let $T$  be a Riesz operator which is not power compact, \cite[Theorem 4.2]{BS} gives an example of such operator.  Then  $ T \in \Phi_{\mathcal{P}B}(X)$ but $ T \notin \Phi_{\mathcal{W}B}(X) $   and so  $\Phi_{\mathcal{W}B}(X) \subsetneq \Phi_{\mathcal{P}B}(X).$

Moreover let $\Phi_{\mathbb{P}}(X)= \{ T \in L(X) \mid \exists P, \,\,\,  P^2=P  \,\, and \,\, T \, is \,\, $P-$ Fredholm \}$ be the class of all $P$-Fredholm operators, when $P$ varies in the set of all idempotents of $L(X).$ Then:  $$ \Phi_{\mathbb{P}}(X)=  \bigcup\limits_{ \{ P\in L(X) \mid P^2=P\}}\Phi_P(X) $$ and   $\Phi_{\mathbb{P}}(X)$  contains strictly  the  class $\Phi_{\mathcal{P}B}(X)$ as shown by the following example.

\begin{example} Let $T: l^2(
\mathbb{N}) \rightarrow l^2(
\mathbb{N}) $  be the forward shift
operator defined by:
$$T(x_1,x_2,x_3,...) = (0, x_1,x_2,x_3,...).$$
Then from \cite[Problem 7.5.3]{AA},   the essential spectrum of $T$ coincides with the unit circle, that is,
$\sigma_e (T)=\{ z \in \mathbb{C}\mid \,\,  \mid z\mid =1 \}.$

Let $ S=T-I,$  then $0 \in \sigma_e (S)$ and $ 0$ is not isolated in $\sigma_e (S).$ Thus $S$ is not a pseudo B-Fredholm operator but $S+I$ is a Fredholm operator and $I$ commutes with $S.$ Hence $ S \in \Phi_I(X)$  and so $ S \in \Phi_{\mathbb{P}}(X) $    but $ S \notin   \Phi_{\mathcal{P}B}(X).$

\end{example}

Thus we have the following strict inclusions between these different  classes :

$$ \Phi(X) \subsetneq \Phi_B(X) \subsetneq \Phi_{\mathcal{W}B}(X) \subsetneq \Phi_{\mathcal{P}B}(X)\subsetneq \Phi_\mathbb{P}(X).$$

Similar inclusion relations can be obtained by considering the left and right versions of the previous classes.

\eremark

\bet \label{left-weak} Let $T \in L(X).$   Then the following properties
are equivalent:

\begin{enumerate}

\item $T$  is  a left weak semi-B-Fredholm operator.

\item  $T$ is a compact perturbation of a power compact-left semi-Fredholm operator.

\end{enumerate}

\eet

\bp  $1)\Rightarrow 2)$  Assume that $T$  is  a left weak semi-B-Fredholm operator. So
$\Pi(T)$  is left Drazin invertible in  $\mathcal{C}(X)$ and  there exist an
idempotent $p= \Pi(P) \in  \mathcal{C}(X)$ such that:

\begin{itemize}
  \item  $p \Pi(T) = \Pi(T) p$
  \item $p \Pi(T)$ is nilpotent in $\mathcal{C}(X)$
  \item  There exists  $S \in L(X)$  such that $p \Pi(S) = \Pi(S) p$ and  $ \Pi(S)( \Pi(T) + p)) = \Pi(I)$

\end{itemize}

As $p \Pi(T)$ is nilpotent in $\mathcal{C}(X),$  then $PT $  is a power compact operator. Let $ X= X_1 \oplus X_2$ be
the decomposition of $X$ associated to $P,$ that's  $ X_1= R(P)$ and $ X_2= N(P).$
 Since $\Pi(T)$ and $\Pi(P)$ commutes, we have that $\Pi(PTP)=\Pi(TP)$ and $\Pi((I-P)T(I-P))=\Pi(T(I-P))$. It follows that $PTP$
 is power compact, $TP=PTP+K_1$, $T(I-P)=(I-P)T(I-P)+K_2$, where $K_1,K_2\in K(X)$, and so:
  \begin{eqnarray*}
    T &=& TP+ T(I-P)=PTP+(I-P)T(I-P)+K,
  \end{eqnarray*}
  \noindent where $K=K_1+K_2\in K(X)$.
    We have that $(X_1, X_2)\in Red(PTP)$ and  $(X_1, X_2)\in Red((I-P)T(I-P))$, hence:
  \begin{equation*}
    PTP = (PTP)_{\mid X_1}\oplus (PTP)_{\mid X_2}=(PTP)_{\mid X_1}\oplus 0,
  \end{equation*}
  and

     $  (I-P)T(I-P) =((I-P)T(I-P)) _{\mid X_1}\oplus ((I-P)T(I-P))_{\mid X_2}=0\oplus ((I-P)T(I-P))_{\mid X_2}.$
Therefore:

\begin{equation}\label{first3}
    T=(PTP)_{\mid X_1}\oplus ((I-P)T(I-P))_{\mid X_2}+K.
  \end{equation}

  It's easily seen that $(PTP)_{\mid X_1}$ is power compact. Moreover we prove that $((I-P)T(I-P))_{\mid X_2}$ is a left  semi-Fredholm operator.
Since $\Pi(S) ( \Pi(T) + p) = \Pi(I),$ we have: $$ (I-P)S(I-P) (I-P) (T+P)(I-P)= I-P  +F_1,$$  where   $F_1$ is compact  such that $ (X_1, X_2) \in Red(F_1).$
As $(I-P)_{\mid R(I-P)}$ is the identity on $X_2$,  it follows   that
$((I-P)T(I-P))_{\mid X_2}=((I-P)(T+P)(I-P))_{\mid X_2}$ is a left  semi-Fredholm operator.
 According to \eqref{first3}, we see that  $T$ is a compact perturbation of
a power compact-left semi-Fredholm operator.\\

$ 2)\Rightarrow 1)$ Assume now that $ X= X_1 \oplus X_2$  is the direct sum of the two closed subspaces $X_1, X_2$  of $X$ and assume that $T=T_1\oplus T_2+ K $ where $T_1 \in L(X_1) $ is a power compact operator, $T_2 \in L(X_2)$ is a left semi-Fredholm operator and $K\in K(X)$. Let $P$ be the projection of $X$ on $X_1$ in parallel to $X_2.$
Then $ PT= (T_1 \oplus 0) +PK$ and $TP= (T_1 \oplus 0) +KP.$  Thus $ PT-TP= PK-KP$ and $PK-KP$  is compact. Let $p= \Pi(P),$  then   $p$ commutes with $ \Pi(T).$
Let also $I_1= {P_|}_{ X_1}$  and $I_2= {(I-P)_|}_{X_2}.$
As $T_2$ is a left semi-Fredholm operator in $L(X_2),$  there exists $ S_2 \in  L(X_2)$  such that $ S_2 T_2- I_2$ is compact, and $I_1  \oplus S_2$ commutes with $P$   because $(I_1  \oplus S_2)P= P(I_1  \oplus S_2)= I_1  \oplus 0. $

We have $ T_1\oplus T_2+P = ( T_1  \oplus 0) +  P  + ( 0 \oplus T_2) = P[ (T_1 \oplus 0)+ I)] P  + (I-P) ( I_1 \oplus  T_2 ) (I-P).$  As $T_1$ is power compact, then $T_1 \oplus 0$ is also power compact and $\Pi( (T_1 \oplus 0)+ I))=  \Pi( (T_1 \oplus 0)) + \Pi(I)$ is invertible in $\mathcal{C}(X).$ Let $ \Pi(S_1) $ be its inverse, where $ S_1 \in L(X).$
We observe that $\Pi( I_1\oplus  T_2) $ is left invertible in  $\mathcal{C}(X)$ having  $ \Pi( I_1 \oplus S_2)$ as a left inverse.
Then:  $$ \Pi(PS_1P + (I-P)( I_1 \oplus S_2) (I-P)) \Pi ( T+P)=\Pi(PS_1P + (I-P)( I_1 \oplus S_2) (I-P)) \Pi ( T_1\oplus T_2+P)$$ $$=\Pi(PS_1P + (I-P)( I_1 \oplus S_2) (I-P)) \Pi ( P (T_1 \oplus 0)+ I) P  + (I-P) ( I_1 \oplus  T_2 ) (I-P))=\Pi(I). $$ It is easily seen that  $\Pi(PS_1P + (I-P)( I_1 \oplus S_2) (I-P))$ commutes with $\Pi(P).$

\noindent Moreover we have $\Pi(PT)= \Pi( T_1 \oplus 0)$   is  nilpotent in $\mathcal{C}(X),$   because $ T_1 \oplus 0$ is a power compact operator.  Thus $ \Pi(T)$ is  left Drazin in $\mathcal{C}(X)$ and  $T$  is  a left weak semi-B-Fredholm operator. \ep

\bigskip

Following the same method as in Theorem \ref{left-weak}, we obtain the following characterization of right pseudo semi-B-Fredholm operators. Due to the similarity of its proof with that of
Theorem \ref{left-weak}, we give it without proof.

\bet \label{right-weak} Let $T \in L(X).$   Then the following properties
are equivalent:

\begin{enumerate}

\item $T$  is  a right weak semi-B-Fredholm operator.

\item  $T$ is a compact perturbation of a power compact-right semi-Fredholm operator.

\end{enumerate}

\eet

Combining Theorem\ref{left-weak} and Theorem\ref{right-weak}, we obtain the following natural characterization of pseudo B-Fredholm operators.

\bet \label{first} Let $T \in L(X).$   Then the following properties
are equivalent:

\begin{enumerate}

\item $T$  is  a  weak B-Fredholm operator.

\item $T$  is  a left and right weak semi-B-Fredholm operator.

\end{enumerate}

\eet

\bp  $1)\Rightarrow 2)$   Assume that  $T$  is  a  weak B-Fredholm operator. So $ \Pi(T)$ is  Drazin invertible in
 $\mathcal{C}(X)$ and from \cite[Proposition 4.9 ]{KP}, there exists an idempotent $p \in \mathcal{C}(X) $ such that $\Pi(T)$
is  $p$-invertible  and
$\Pi(T)p$ is nilpotent. So  $ \Pi(T)$ is  left  Drazin invertible  and right  Drazin invertible  in  $\mathcal{C}(X)$ and $T$  is  a left and right weak semi-B-Fredholm operator.

$2)\Rightarrow 1)$ Conversely assume now that $\Pi(T)$  is left and  right  Drazin invertible in  $\mathcal{C}(X),$  so  $\Pi(T)$  is left and  right generalized Drazin invertible in  $\mathcal{C}(X).$ From Theorem \ref{first equivalence}, it follows that  $\Pi(T)$ is generalized Drazin invertible in  $\mathcal{C}(X).$ Then from \cite[Theorem 3.1]{Koliha} there exists an idempotent $p \in \mathcal{C}(X),$ commuting with $ \Pi(T)$   such that $\Pi(T)p$ is quasi-nilpotent and $ \Pi(T)+p$ is invertible. From \cite[Theorem 3.1]{Koliha}, we know that $p$ lies in the closed subalgebra of $ \mathcal{C}(X)$   generated by $\Pi(T).$

Now as  $ \Pi(T)$ is  left  Drazin invertible in $\mathcal{C}(X),$ there exists an idempotent $ p_0 \in  \mathcal{C}(X),$ commuting with $\Pi(T)$ such that $\Pi(T)p_0$ is nilpotent and $\Pi(T) + p_0$ is let invertible. Assume that $ (\Pi(T)p_0)^n= (\Pi(T))^n p_0=  0,$  for an integer $n>0.$ Since $p_0$ commutes with $\Pi(T),$  it commutes   also with $p.$  We have: $$ (\pi(T))^n p =  [(\pi(T))^n p_0 +  (\pi(T))^n (e-p_0)]p=[(\pi(T))^n (e-p_0)]p, $$  where $ e= \Pi(I)$ is the identity element of $\mathcal{C}(X).$ From \cite[Theorem 2.15]{P50}, it follows that  $(\pi(T))^n (e-p_0)$ is left invertible in the closed subalgebra
  $(e-p_0)\mathcal{C}(X) (e-p_0)$   of $\mathcal{C}(X).$ So there exists $ b \in (e-p_0)\mathcal{C}(X) (e-p_0)$ such that $ b(\pi(T))^n (e-p_0))= e-p_0.$  Hence: $$ b^m((\pi(T))^n )^m (e-p_0)p = (e-p_0)p\,\, \text{and} \,\, b^m((\pi(T))^n p)^m(e-p_0)  = p(e-p_0), $$  for all integer $m > 0.$
As $ p$ commutes with $p_0,$  then  $p(e-p_0)$ is an idempotent.  If $  p(e-p_0) \neq 0,$  then $\norm{(e-p_0)p} \geq 1.$   So $  1 \leq \norm{b} \norm{(\pi(T)^n p )^m}^{\frac{1}{m}} {\norm{e-p_0}}^{\frac{1}{m}} $  for all integer $m > 0.$  However this is impossible since $ (\pi(T))^n p$ is quasinilpotent.
  So necessarily  $ (e-p_0)p = 0 $ and   $((\pi(T))^n p)= 0.$  Thus $ \pi(T)p$ is nilpotent and $ \Pi(T)$ is Drazin invertible in   $\mathcal{C}(X).$ Therefore  $T$  is  a  weak B-Fredholm operator.
\ep

  \noindent  Alaa Hamdan\\
\noindent   Dubai \\
\noindent   United Arab Emirates \\
\noindent     aa.hamdan@outlook.com \\

\noindent  Mohammed Berkani\\
\noindent Honorary member of LIAB \\
\noindent Department of Mathematics\\
\noindent Science Faculty of Oujda \\
\noindent Mohammed I University\\
\noindent Morocco\\
\noindent berkanimo@aim.com

\end{document}